\documentclass[12pt,english]{article}
\usepackage[latin9]{inputenc}
\usepackage{geometry}
\usepackage{amssymb}
\usepackage{wasysym}

\makeatletter
\date{}

\makeatother

\usepackage{babel}
\begin{document}

\title{Induced Saturation of $P_{6}$}

\author{Eero R{\"a}ty\thanks{Centre for Mathematical Sciences, Wilberforce Road, Cambridge CB3
0WB, UK, epjr2@cam.ac.uk}}
\maketitle
\begin{abstract}
A graph $G$ is called $H$-induced-saturated if $G$ does not contain
an induced copy of $H$, but removing any edge from $G$ creates an
induced copy of $H$ and adding any edge of $G^{c}$ to $G$ creates
an induced copy of $H$. Martin and Smith showed that there does not
exist a $P_{4}$-induced-saturated graph, where $P_{4}$ is the path
on 4 vertices. Axenovich and Csik{\'o}s studied related questions, and
asked if there exists a $P_{n}$-induced-saturated graph for any $n\geq5$.
Our aim in this short note is to show that there exists a $P_{6}$-induced-saturated
graph.
\end{abstract}

\section{Introduction}

A graph $G$ is said to be \textit{$H$-saturated} if $G$ does not
contain a copy of $H$, but adding any new edge to $G$ gives a copy
of $H$. It is clear that for any non-empty $H$ there exists such
$G$ with a given number of vertices.

The notion of saturation can be generalised to induced subgraphs in
the following way. A graph $G$ is said to be \textit{$H$-induced-saturated}
if $G$ does not contain an induced copy of $H$, but removing any
edge from $G$ creates an induced copy of $H$ and adding any edge
of $G^{c}$ to $G$ creates an induced copy of $H$. 

It is not clear anymore whether for given $H$ there exists a graph
$G$ which is $H$-induced-saturated. In fact, it was proved by Martin
and Smith \cite{key-2} that for $H=P_{4}$ there is no such $G$
satisfying the property, where $P_{n}$ denotes the path with $n$
vertices. For convenience, we say that $H$ is \textit{induced-saturated}
if there exists some $G$ which is $H$-induced-saturated. 

It is natural to ask what happens when $H=P_{n}$ for other values
of $n$. The cases $H=P_{2}$ and $H=P_{3}$ are trivial, as one can
take $G$ to be empty graph and $K_{m}$ respectively. Axenovich and
Csik{\'o}s \cite{key-3} investigated several related questions by giving
examples of families of graphs that are induced-saturated, and also
asked if the graphs $H=P_{n}$ are induced-saturated for $n\geq5$.
The aim of this note is to provide an example which shows that $P_{6}$
is induced-saturated. \\

\section{The construction}

\textbf{Theorem 1. }There exists $G$ which is $P_{6}$-induced-saturated. 
\\

\textbf{Proof.} Let $\mathbb{F}=\mathbb{F}_{2}\left(\alpha\right)/\left(\alpha^{4}+\alpha+1\right)$
be the finite field of order 16, and note that $\alpha$ is a generator
for the multiplicative group $\mathbb{F}^{\times}$. Let $S$ be the
set of non-zero cubes in $\mathbb{F}$, i.e. $S=\left\{ 1,\,\alpha^{3},\,\alpha^{2}+\alpha^{3},\,\alpha+\alpha^{3},\,1+\alpha+\alpha^{2}+\alpha^{3}\right\} $.
Define a graph $G$ whose vertex set is $\mathbb{F}$ and whose edges
are given by $xy\in E\left(G\right)$ if and only $x-y\in S$. 

First of all, note that $x\rightarrow\alpha^{3i}x+\beta$ is an automorphism
of $G$ for any $i\in\left\{ 0,\dots,4\right\} $ and $\beta\in\mathbb{F}$.
Thus for any two edges $e,\,f\in E\left(G\right)$ there exists an
automorphism $\theta$ such that $\theta\left(e\right)=f$, i.e. the
group of automorphisms acts transitively on the edges of $G$. 

We will now check that $G$ satisfies all the required properties.
\\

\textbf{Claim 1. }$G$ does not contain an induced copy of $P_{6}$.
\\

\textbf{Proof.} Suppose $G$ contains copy of $P_{6}$. Since the
group of automorphisms acts transitively on $E\left(G\right)$, we
may assume that the first two vertices are $0$ and $1$. Let $T=\Gamma\left(0\right)^{c}\cap\Gamma\left(1\right)^{c}$.
Then the last three vertices on the path are elements of $T$. It
is easy to verify that 
\[
T=\left\{ \alpha,\,1+\alpha,\,\alpha^{2},\,1+\alpha^{2},\,\alpha+\alpha^{2},\,1+\alpha+\alpha^{2}\right\} .
\]
Hence $G\left[T\right]$ is union of three disjoint edges corresponding
to the pairs $\left\{ \alpha,1+\alpha\right\} $, $\left\{ \alpha^{2},1+\alpha^{2}\right\} $
and $\left\{ \alpha+\alpha^{2},\,1+\alpha+\alpha^{2}\right\} $. Thus
$G\left[T\right]$ does not contain a path with three vertices, so
$G$ cannot contain an induced copy of $P_{6}$. $\square$\\

\textbf{Claim 2.} Adding a new edge to $G$ gives an induced copy
of $P_{6}$.\\

\textbf{Proof. }Note that for any edge $e=\left\{ x,y\right\} $ in
the complement of $G$ there is an automorphism of $G$ which maps
$e$ to either $\left\{ 0,\alpha^{10}\right\} $ or $\left\{ 0,\alpha^{14}\right\} $
(corresponding to the cases $x-y=\alpha^{3i+1}$ and $x-y=\alpha^{3i+2}$
respectively). Hence it suffices to consider just the cases $e=\left\{ 0,\alpha^{10}\right\} $
and $e=\left\{ 0,\alpha^{14}\right\} $. 

Consider the particular elements $x_{1}=0$, $x_{2}=\alpha+\alpha^{3}$,
$x_{3}=\alpha$, $x_{4}=\alpha+\alpha^{2}+\alpha^{3}$, $x_{5}=\alpha^{2}$,
$x_{6}=\alpha^{10}=1+\alpha+\alpha^{2}$ and $x_{7}=\alpha^{14}=1+\alpha^{3}$,
and let $R=\left\{ x_{1},\dots,x_{7}\right\} $. It is easy to verify
that $G\left[R\right]$ is union of induced $P_{5}$ with vertices
$x_{1},\,x_{2},\,x_{3},\,x_{4},\,x_{5}$ (in this order) together
with two isolated vertices $x_{6}$ and $x_{7}$. Hence addition of
either of $\left\{ 0,\alpha^{10}\right\} $ or $\left\{ 0,\alpha^{14}\right\} $
gives an induced copy of $P_{6}$ in $G$. $\Square$\\

\textbf{Claim 3. }Removal of any edge from $G$ gives an induced copy
of $P_{6}$.\\

\textbf{Proof.} For any $e\in E\left(G\right)$ there is an automorphism
of $G$ mapping $e$ to the edge $\left\{ 0,1\right\} $. Hence it
suffices to check that $G$ with edge $f=\left\{ 0,1\right\} $ removed
contains an induced copy of $P_{6}$. It is easy to check that $v_{1}=\alpha+\alpha^{2}+\alpha^{3}$,
$v_{2}=1$, $v_{3}=1+\alpha^{3}$, $v_{4}=\alpha^{3}$, $v_{5}=0$
and $v_{6}=\alpha^{2}+\alpha^{3}$ forms an induced copy of $P_{6}$.
$\square$\\

From the claims above it follows that $G$ is $P_{6}$-induced-saturated
graph. $\square$\\

Unfortunately the construction gives no idea on what happens for other
values of $n$. Note that the conditions 'removing any edge from $G$
gives an induced $P_{n}$' and 'adding any missing edge to $G$ gives
an induced $P_{n}$' can be together written as 'there exists $f\,:\,V\left(G\right)^{\left(2\right)}\rightarrow V\left(G\right)^{\left(n-2\right)}$
such that for any $e\in V\left(G\right)^{\left(2\right)}$, we have
$e\cap f\left(e\right)=\emptyset$ and $G_{e}\left[e\cup f\left(e\right)\right]$
is isomorphic to $P_{n}$', where $G_{e}$ is the graph obtained from
$G$ by removing or adding $e$ depending on whether $e\in E\left(G\right)$
or $e\not\in E\left(G\right)$ respectively. If $m=\left|V\left(G\right)\right|$,
then $f$ can be viewed as a function $f\,:\,\left[m\right]^{\left(2\right)}\rightarrow\left[m\right]^{\left(n-2\right)}$.
Even though the condition that $G$ does not contain induced copy
of $P_{n}$ is not yet considered, this suggests that the other conditions
should be easier to satisfy when $n$ increases. Hence it might be
reasonable to suppose that for $n>6$ such a construction exists.

\end{document}